\documentclass[11pt, reqno]{amsart} \textwidth=14.9cm \oddsidemargin=1cm \evensidemargin=1cm
\usepackage{amsmath,amssymb,amscd,mathrsfs,stmaryrd,wasysym}
\usepackage{amsmath}
\usepackage{amsxtra}
\usepackage{amscd}
\usepackage{amsthm}
\usepackage{amsfonts}
\usepackage{amssymb}
\usepackage{eucal}
\usepackage{color}
\usepackage{graphicx}%
\usepackage{bm}
\newtheorem{theorem}{Theorem}[section]
\newtheorem{lemma}[theorem]{Lemma}
\newtheorem{proposition}[theorem]{Proposition}

\theoremstyle{definition}

\newtheorem{example}[theorem]{Example}

\definecolor{A}{rgb}{.75,1,.75}

\numberwithin{equation}{section}

\begin{document}

\title[Fusion procedure]{Fusion procedure for Yokonuma-Hecke algebras}
\author[Weideng Cui]{Weideng Cui}
\address{School of Mathematics, Shandong University, Jinan, Shandong 250100, P.R. China.}
\email{cwdeng@amss.ac.cn}

\begin{abstract}
Inspired by the work [PA1], in this note, we prove that a complete set of pairwise orthogonal primitive idempotents of Yokonuma-Hecke algebras can be obtained by consecutive evaluations of a certain rational function.
\end{abstract}

%We establish an equivalence between a module category of the affine (resp. cyclotomic) Yokonuma-Hecke algebra (resp. $Y_{r,n}^{\lambda}(q)$) and its suitable counterpart for a direct sum of tensor products of affine Hecke algebras of type $A$ (resp. cyclotomic Hecke algebras). We then develop several applications of this result. The simple modules of affine Yokonuma-Hecke algebras and of their associated cyclotomic Yokonuma-Hecke algebras are classified over an algebraically closed field of characteristic $p=0$ or $(p,r)=1.$ The modular branching rules for these algebras are obtained, and they are further identified with crystal graphs of integrable modules for affine Lie algebras of type $A.$\end{abstract}

%\thanks{{Keywords}: Affine Yokonuma-Hecke algebras; Cyclotomic Yokonuma-Hecke algebras; Modular branching rules; Crystal graphs} \large

\maketitle
\medskip
\section{Introduction}
\subsection{}
Yokonuma-Hecke algebras were introduced by Yokonuma \cite{Yo} as a centralizer algebra associated to the permutation representation of a finite Chevalley group $G$ with respect to a maximal unipotent subgroup of $G$. The Yokonuma-Hecke algebra $Y_{r,n}(q)$ (of type $A$) is a quotient of the group algebra of the modular framed braid group $(\mathbb{Z}/r\mathbb{Z})\wr B_{n},$ where $B_{n}$ is the braid group on $n$ strands (of type $A$). By the presentation given by Juyumaya and Kannan \cite{Ju1, Ju2, JuK}, the Yokonuma-Hecke algebra $Y_{r,n}(q)$ can also be regraded as a deformation of the group algebra of the complex reflection group $G(r,1,n),$ which is isomorphic to the wreath product $(\mathbb{Z}/r\mathbb{Z})\wr \mathfrak{S}_{n}$, where $\mathfrak{S}_n$ is the symmetric group. It is well-known that there exists another deformation of the group algebra of $G(r,1,n),$ namely the Ariki-Koike algebra \cite{AK}. The Yokonuma-Hecke algebra $Y_{r,n}(q)$ is quite different from the Ariki-Koike algebra.  For example, the Iwahori-Hecke algebra of type $A$ is canonically a subalgebra of the Ariki-Koike algebra, whereas it is an obvious quotient of $Y_{r,n}(q),$ but not an obvious subalgebra of it.

Recently, by generalizing the approach of Okounkov-Vershik \cite{OV} on the representation theory of $\mathfrak{S}_n$, Chlouveraki and Poulain d'Andecy \cite{ChPA1} introduced the notion of affine Yokonuma-Hecke algebra $\widehat{Y}_{r,n}(q)$ and gave explicit formulas for all irreducible representations of $Y_{r,n}(q)$ over $\mathbb{C}(q)$, and obtained a semisimplicity criterion for it. In their subsequent paper [ChPA2], they studied the representation theory of the affine Yokonuma-Hecke algebra $\widehat{Y}_{r,n}(q)$ and the cyclotomic Yokonuma-Hecke algebra $Y_{r,n}^{d}(q)$. In particular, they gave the classification of irreducible representations of $Y_{r,n}^{d}(q)$ in the generic semisimple case. In \cite{CW}, we gave the classification of the simple $\widehat{Y}_{r,n}(q)$-modules as well as the classification of the simple modules of the cyclotomic Yokonuma-Hecke algebras over an algebraically closed field $\mathbb{K}$ of characteristic $p$ such that $p$ does not divide $r.$ In the past several years, the study of affine and cyclotomic Yokonuma-Hecke algebras has made substantial progress; see [ChPA1, ChPA2, ChS, C1, C2, CW, ER, JaPA, Lu, PA2, Ro].

\subsection{}
Jucys [Juc] claimed that the primitive idempotents of $\mathfrak{S}_n$ indexed by standard Young tableaux can be obtained by taking a certain limiting process on a rational function, which is now commonly referred to as the fusion procedure. It has been developed in the situation of Hecke algebras [Ch], see also [Na1-3]. Molev [Mo] has proposed an alternative approach of the fusion procedure for the symmetric group, which relies on the existence of a maximal commutative subalgebra generated by the Jucys-Murphy elements. Here the idempotents are derived by consecutive evaluations of a certain rational function. The simple version of the fusion procedure has been generalized to the Hecke algebras of type $A$ [IMO], to the Brauer algebras [IM, IMOg1], to the Birman-Murakami-Wenzl algebras [IMOg2], to the complex reflection groups of type $G(d,1,n)$ [OgPA1], to the Ariki-Koike algebras [OgPA2], to the wreath products of finite groups by the symmetric group [PA1], to the degenerate cyclotomic Hecke algebras [ZL].

\subsection{}
Inspired by [PA1], in this note we prove that a complete set of pairwise orthogonal primitive idempotents for Yokonuma-Hecke algebras can be constructed by consecutive evaluations of a certain rational function.

Assume that $Y_{r,n}$ is the split semisimple Yokonuma-Hecke algebra defined over a field $\mathbb{K}$. We denote by $\mathcal{P}_{r,n}$ the set of all $r$-partitions of $n$ and $\text{STab}_{r}(n)$ the set of all standard $r$-tableaux of size $n$ (of any shape). Then a complete set of pairwise orthogonal primitive idempotents of $Y_{r,n}$ is parameterized by $\text{STab}_{r}(n)$. For each $\bm{\lambda}\in \mathcal{P}_{r,n}$ and a standard $r$-tableau $\mathcal{T}$ of shape $\bm{\lambda},$ let $E_{\mathcal{T}}$ be the primitive idempotent of $Y_{r,n}$ corresponding to $\mathcal{T}.$

We now state the main result of this paper.\vskip0.25cm
\noindent\textbf{Theorem.} The idempotent $E_{\mathcal{T}}$ of $Y_{r,n}$ can be obtained by the following consecutive evaluations \[E_{\mathcal{T}}=\frac{1}{\text{F}_{\bm{\lambda}}^{T}\text{F}_{\bm{\lambda}}}\Phi(u_1,\ldots,u_{n},v_1,\ldots,v_{n})
\Big|_{v_{1}=\zeta_{\text{p}_{1}}}\cdots\Big|_{v_{n}=\zeta_{\text{p}_{n}}}\Big|_{u_{1}=\text{c}_1}\cdots\Big|_{u_{n}=\text{c}_{n}}.\]

This paper is organized as follows. In Section 2, we introduce the Jucys-Murphy elements of $Y_{r, n}$ and recall some combinatorial notions. In Section 3, we recall the construction of the primitive idempotents $E_{\mathcal{T}}$ of $Y_{r, n}$ following [ChPA1], and then establish the fusion formula for $E_{\mathcal{T}}.$

\section{Preliminaries}
In this section, we first introduce the Jucys-Murphy elements of the Yokonuma-Hecke algebra, and then recall some combinatorial notions.

\subsection{Jucys-Murphy elements}
Let $r, n\in \mathbb{N},$ $r\geq1,$ and let $\zeta=e^{2\pi i/r}.$ Let $q$ be an indeterminate.

Let $\mathcal{R}=\mathbb{Z}[\frac{1}{r}][q,q^{-1}].$ The Yokonuma-Hecke algebra $Y_{r,n}=Y_{r,n}(q)$ is an $\mathcal{R}$-associative algebra generated by the elements $t_{1},\ldots,t_{n},g_{1},\ldots,g_{n-1}$ satisfying the following relations:
\begin{equation}\label{rel-def-Y1}\begin{array}{rclcl}
g_ig_j\hspace*{-7pt}&=&\hspace*{-7pt}g_jg_i && \mbox{for all $i,j=1,\ldots,n-1$ such that $\vert i-j\vert \geq 2$;}\\[0.1em]
g_ig_{i+1}g_i\hspace*{-7pt}&=&\hspace*{-7pt}g_{i+1}g_ig_{i+1} && \mbox{for all $i=1,\ldots,n-2$;}\\[0.1em]
t_it_j\hspace*{-7pt}&=&\hspace*{-7pt}t_jt_i &&  \mbox{for all $i,j=1,\ldots,n$;}\\[0.1em]
g_it_j\hspace*{-7pt}&=&\hspace*{-7pt}t_{s_i(j)}g_i && \mbox{for all $i=1,\ldots,n-1$ and $j=1,\ldots,n$;}\\[0.1em]
t_i^r\hspace*{-7pt}&=&\hspace*{-7pt}1 && \mbox{for all $i=1,\ldots,n$;}\\[0.2em]
g_{i}^{2}\hspace*{-7pt}&=&\hspace*{-7pt}1+(q-q^{-1})e_{i}g_{i} && \mbox{for all $i=1,\ldots,n-1$,}
\end{array}
\end{equation}
where $s_{i}$ is the transposition $(i,i+1)$, and for each $1\leq i\leq n-1$,
$$e_{i} :=\frac{1}{r}\sum\limits_{s=0}^{r-1}t_{i}^{s}t_{i+1}^{-s}.$$

Note that the elements $e_{i}$ are idempotents in $Y_{r,n}.$ The elements $g_{i}$ are invertible, with the inverse given by
\begin{equation}\label{inverse}
g_{i}^{-1}=g_{i}-(q-q^{-1})e_{i}\quad\mbox{for~all}~i=1,\ldots,n-1.
\end{equation}

Let $w\in \mathfrak{S}_{n},$ and let $w=s_{i_1}\cdots s_{i_{r}}$ be a reduced expression of $w.$ By Matsumoto's lemma, the element $g_{w} :=g_{i_1}g_{i_2}\cdots g_{i_{r}}$ does not depend on the choice of the reduced expression of $w,$ that is, it is well-defined.

Let $i, k\in \{1,2,\ldots,n\}$ and set \begin{equation}\label{idempotents}e_{i,k} :=\frac{1}{r}\sum\limits_{s=0}^{r-1}t_{i}^{s}t_{k}^{-s}.\end{equation} Note that $e_{i,k}^{2}=e_{i,k}=e_{k,i},$ and that $e_{i,i+1}=e_{i}.$ It can be easily checked that
\begin{equation}\label{relations}\begin{array}{rclcl}
t_{i}e_{j,k}\hspace*{-7pt}&=&\hspace*{-7pt}e_{j,k}t_{i} && \mbox{for all $i,j,k=1,\ldots,n$,}\\[0.1em]
e_{i,j}e_{k,l}\hspace*{-7pt}&=&\hspace*{-7pt}e_{k,l}e_{i,j} &&  \mbox{for all $i,j,k,l=1,\ldots,n$,}\\[0.1em]
e_ie_{k,l}\hspace*{-7pt}&=&\hspace*{-7pt}e_{s_i(k), s_i(l)}e_i && \mbox{for all $i=1,\ldots,n-1$ and $k,l=1,\ldots,n$,}\\[0.1em]
e_{j,k}g_{i}\hspace*{-7pt}&=&\hspace*{-7pt}g_{i}e_{s_{i}(j),s_{i}(k)} && \mbox{for all $i=1,\ldots,n-1$ and $j,k=1,\ldots,n$.}
\end{array}
\end{equation}
In particular, we have $e_{i}g_{i}=g_{i}e_{i}$ for all $i=1,2,\ldots,n-1.$

We define inductively the following elements in $Y_{r,n}$:
\begin{equation}\label{JM-elements}
J_{1} :=1\text{ and } J_{i+1} :=g_{i}J_ig_i\quad\mbox{for}~i=1,\ldots,n-1.
\end{equation}
By [ChPA1, Corollary 1] we have, for any $1\leq i\leq n-1$,
\begin{equation}\label{giXj}
g_{i}J_{j}=J_{j}g_{i}\quad\mathrm{for}~j=1,2,\ldots,n~\mathrm{such~that}~j\neq i, i+1.
\end{equation}

Let $\mathbb{K}$ be an algebraically closed field of characteristic $p\geq 0$ such that $p$ does not divide $r$. In the rest of this paper, we shall work with a specialised split semisimple Yokonuma-Hecke algebra $Y_{r, n}$ defined over $\mathbb{K},$ that is, $q\in \mathbb{K}^{*}$ satisfies the following condition:
\begin{equation}\label{semisimple}
\prod_{k=1}^{n}(1+q^{2}+\cdots+q^{2(k-2)}+q^{2(k-1)})\neq 0.
\end{equation}

The elements $J_1,\ldots,J_n,$ together with the elements $t_1,\ldots,t_n,$ are called the Jucy-Murphy elements of $Y_{r,n}$, which generate a maximal commutative subalgebra of $Y_{r,n}.$

\subsection{$r$-partitions}
$\lambda=(\lambda_{1},\ldots,\lambda_{k})$ is called a partition of $n$ if it is a finite sequence of non-increasing nonnegative integers whose sum is $n.$ We write $\lambda\vdash n$ if $\lambda$ is a partition of $n,$ and we set $|\lambda| :=n$. We associate a Young diagram to a partition $\lambda$, which is the set $$[\lambda] :=\{(i,j)\:|\:i\geq 1~\mathrm{and}~1\leq j\leq \lambda_{i}\}.$$ We will regard $[\lambda]$ as a left-justified array of rows containing $\lambda_{j}$ nodes in the $j$-th row for $j=1,\ldots,k.$ We write $\theta=(x,y)$ for the node in row $x$ and column $y.$

For a partition $\lambda,$ a node $\theta\in [\lambda]$ is called removable from $\lambda$ if the set of nodes obtained from $[\lambda]$ by removing $\theta$ is still a Young diagram of a partition. A node $\theta'\notin [\lambda]$ is called addable to $\lambda$ if the set of nodes obtained from $[\lambda]$ by adding $\theta'$ is still a Young diagram of a partition. The conjugate of $\lambda$ is the partition $\lambda'=(\lambda'_1,\ldots,\lambda'_l)$, which is given by $$\lambda'_j :=\sharp\{i\:|\:1\leq i\leq k \text{ such that }\lambda_i\geq j\}.$$

An $r$-partition of $n$ is an ordered $r$-tuple $\bm{\lambda}=(\lambda^{(1)},\lambda^{(2)},\ldots,\lambda^{(r)})$ of partitions $\lambda^{(k)}$ such that $\sum_{k=1}^{r}|\lambda^{(k)}|=n.$ We denote by $\mathcal{P}_{r,n}$ the set of $r$-partitions of $n.$ A pair $\bm{\theta}=(\theta, k)$ consisting of a node $\theta$ and an integer $k\in \{1,\ldots,r\}$ is called an $r$-node, and the integer $k$ is called the position of $\bm{\theta}.$ The Young diagram $[\bm{\lambda}]$ of an $r$-partition $\bm{\lambda}$ is the ordered $r$-tuple of the Young diagram of its components.

Let $\bm{\lambda}=(\lambda^{(1)},\lambda^{(2)},\ldots,\lambda^{(r)})$ be an $r$-partition. An $r$-node $\bm{\theta}=(\theta, k)\in [\bm{\lambda}]$ is called removable from $\bm{\lambda}$ if the node $\theta$ is removable from $\lambda^{(k)}$. An $r$-node $\bm{\theta'}=(\theta', k')\notin [\bm{\lambda}]$ is called removable to $\bm{\lambda}$ if the node $\theta'$ is addable to $\lambda^{(k')}$.

For an $r$-node $\bm{\theta}=((x,y),k)$, we define $\text{cc}(\bm{\theta}) :=y-x,$ $\text{p}(\bm{\theta}) :=k$ and the quantum content $\text{c}(\bm{\theta}) :=q^{2(y-x)}.$

\subsection{Hook length}
Let $\bm{\lambda}=(\lambda^{(1)},\lambda^{(2)},\ldots,\lambda^{(r)})$ be an $r$-partition and $\bm{\theta}=(\theta, k)=((x,y),k)$ an $r$-node of $[\bm{\lambda}].$ We define the hook length $h_{\bm{\lambda}}(\bm{\theta})$ of $\bm{\theta}$ in $\bm{\lambda}$ to be the hook length of the node $\theta$ in $\lambda^{(k)},$ that is,
\begin{equation}\label{hook-length}
h_{\bm{\lambda}}(\bm{\theta}) :=h_{\lambda^{(k)}}(\theta)=\lambda^{(k)}_{x}+\lambda^{(k)'}_{y}-x-y+1.
\end{equation}

Set $\zeta_{k} :=\zeta^{k-1}$ for $1\leq k\leq r,$ and let $S :=\{\zeta_{1},\zeta_{2},\ldots,\zeta_{r}\}.$ For an $r$-partition $\bm{\lambda}=(\lambda^{(1)},\lambda^{(2)},\ldots,\lambda^{(r)})$, we define
\begin{equation}\label{f-lambda-t}
\text{F}_{\bm{\lambda}}^{T} :=\prod_{\bm{\theta}\in \bm{\lambda}} \Big(\prod_{\substack{\xi\in S\\\xi\neq \zeta_{\text{p}(\bm{\theta})}}}(\zeta_{\text{p}(\bm{\theta})}-\xi)\Big),
\end{equation}
and
\begin{equation}\label{f-lambda}
\text{F}_{\bm{\lambda}} :=\prod_{\bm{\theta}\in \bm{\lambda}} \frac{[h_{\bm{\lambda}}(\bm{\theta})]_{q}}{q^{\text{cc}(\bm{\theta})}}=\prod_{k=1}^{r}\prod_{\theta\in [\lambda^{(k)}]}\frac{[h_{\lambda^{(k)}}(\theta)]_{q}}{q^{\text{cc}(\bm{\theta})}},
\end{equation}
where $[a]_{q}=q^{a-1}+q^{a-3}+\cdots+q^{-a+1}$ for $a\in \mathbb{Z}_{\geq 0}.$

\subsection{Standard $r$-tableaux}
Let $\bm{\lambda}$ be an $r$-partition of $n.$ An $r$-tableau of shape $\bm{\lambda}$ is a bijection between the set $\{1,\ldots,n\}$ and the set of $r$-nodes in $[\bm{\lambda}]$, and the number $n$ is called the size of the $r$-tableau. An $r$-tableau is called standard if the numbers increase along any row (from left to right) and down any column (from top to bottom) of each diagram in $[\bm{\lambda}].$ Denote by $\mathrm{Std}(\bm{\lambda})$ the set of standard $r$-tableaux of shape $\bm{\lambda}.$

Let $\bm{\lambda}$ be an $r$-partition of $n$ and $\mathcal{T}$ a standard $r$-tableau of shape $\bm{\lambda}$. We denote by $\text{c}(\mathcal{T}|i)$ and $\text{p}(\mathcal{T}|i)$ the quantum content and the position of the $r$-node containing the integer $i,$ respectively. For brevity, we set
\begin{equation}\label{ci-pi-12345}
\text{c}_i :=\text{c}(\mathcal{T}|i)\quad\text{ and } \quad\text{p}_i :=\text{p}(\mathcal{T}|i)\quad\text{ for }i=1,\ldots,n.
\end{equation}
We then define
\begin{equation}\label{f-TT}
\text{F}_{\mathcal{T}}^{T}(v) :=\prod_{\substack{\xi\in S\\\xi\neq \zeta_{\text{p}_{n}}}}\frac{1}{v-\xi},
\end{equation}
and
\begin{equation}\label{f-T}
\text{F}_{\mathcal{T}}(u) :=\frac{u-\text{c}_n}{u-1}\prod_{i=1}^{n-1}\frac{(u-\text{c}_i)^{2}}{(u-\text{c}_i)^{2}-(q-q^{-1})^{2}u\text{c}_i\delta_{\text{p}_i, \text{p}_n}},
\end{equation}
where $\delta_{\text{p}_i, \text{p}_n}$ is the Kronecker delta.

Let $\bm{\mu}$ be the shape of the standard $r$-tableau obtained from $\mathcal{T}$ by removing the $r$-node containing the number $n$. Recall that for any fixed $r$-th root of unity $\xi$, we have
\begin{equation}
\label{young-module-5-19}\prod_{\xi\neq \alpha\in S}(\xi-\alpha)=r\xi^{-1}.
\end{equation}
Thus, $\text{F}_{\mathcal{T}}^{T}(v)$ is non-singular at $v=\zeta_{\text{p}_{n}}.$ Moreover, from \eqref{f-lambda-t} we have
\begin{equation}\label{f-T-relation}
\text{F}_{\mathcal{T}}^{T}(v)\Big|_{v=\zeta_{\text{p}_{n}}}=\frac{\zeta_{\text{p}_{n}}}{r}=(\text{F}_{\bm{\lambda}}^{T})^{-1}\text{F}_{\bm{\mu}}^{T}.
\end{equation}

The following proposition can be proved in exactly the same way as in [OgPA1, Propositions 3.4 and 4.4].
\begin{proposition}
The rational function $\emph{F}_{\mathcal{T}}(u)$ is non-singular at $u=\emph{c}_n$, and moreover, we have
\begin{equation}\label{f-T-u}
\emph{F}_{\mathcal{T}}(u)\Big|_{u=\emph{c}_n}=\emph{F}_{\bm{\lambda}}^{-1}\emph{F}_{\bm{\mu}}.
\end{equation}
\end{proposition}

\subsection{Baxterized elements}
We define the following rational functions in variables $a,b$ with values in $Y_{r,n}$:
\begin{equation}\label{Baxter-element}
g_{i}(a,b) :=g_{i}+(q-q^{-1})\frac{be_{i}}{a-b}\quad\mbox{for}~i=1,\ldots,n-1.
\end{equation}
The functions $g_{i}(a,b)$ are called Baxterized elements. We then have following lemma.
\begin{lemma}\label{Bax-elements}
The Baxterized elements $g_{i}(a,b)$ satisfy the following relations$:$
\begin{equation}\label{Baxter-element1}
g_{i}(a,b)g_{i+1}(a,c)g_{i}(b,c)=g_{i+1}(b,c)g_{i}(a,c)g_{i+1}(a,b)\quad\mbox{for}~i=1,\ldots,n-1,
\end{equation}
\begin{equation}\label{Baxter-element2}
\hspace*{20pt}g_{i}(a,b)g_{i}(b,a)=1-(q-q^{-1})^{2}\frac{abe_{i}}{(a-b)^{2}}\qquad\mbox{for}~i=1,\ldots,n-1.
\end{equation}
\end{lemma}
\begin{proof}
We first prove \eqref{Baxter-element2}. By \eqref{relations} we have
\begin{align*}
g_{i}(a,b)g_{i}(b,a)&=g_{i}^{2}+(q-q^{-1})\frac{ag_{i}e_i}{b-a}+(q-q^{-1})\frac{be_{i}g_i}{a-b}+(q-q^{-1})^{2}\frac{abe_{i}^{2}}{(a-b)(b-a)}\\
&=1+(q-q^{-1})e_{i}g_{i}-(q-q^{-1})g_{i}e_{i}-(q-q^{-1})^{2}\frac{abe_{i}}{(a-b)^{2}}\\
&=1-(q-q^{-1})^{2}\frac{abe_{i}}{(a-b)^{2}}.
\end{align*}

Next we prove \eqref{Baxter-element1}. By definition, we see that the left-hand side of \eqref{Baxter-element1} is equal to
\begin{align*}
\big(g_{i}+(q-q^{-1})\frac{be_{i}}{a-b}\big)\big(g_{i+1}+(q-q^{-1})\frac{ce_{i+1}}{a-c}\big)\big(g_{i}+(q-q^{-1})\frac{ce_{i}}{b-c}\big).
\end{align*}
By expanding the expression above, we get that it is equal to
\begin{align}\label{equation1}
g_{i}g_{i+1}g_{i}+&(q-q^{-1})\frac{cg_{i}g_{i+1}e_{i}}{b-c}+(q-q^{-1})\frac{cg_{i}e_{i+1}g_{i}}{a-c}\notag\\
+(q&-q^{-1})^{2}\frac{c^{2}g_{i}e_{i+1}e_{i}}{(a-c)(b-c)}+
(q-q^{-1})\frac{be_{i}g_{i+1}g_{i}}{a-b}+(q-q^{-1})^{2}\frac{bce_{i}g_{i+1}e_{i}}{(a-b)(b-c)}\notag\\
&+(q-q^{-1})^{2}\frac{bce_{i}e_{i+1}g_{i}}{(a-b)(a-c)}
+(q-q^{-1})^{3}\frac{bc^{2}e_{i}e_{i+1}e_{i}}{(a-b)(a-c)(b-c)}.
\end{align}

Similarly, we see that the right-hand side of \eqref{Baxter-element1} is equal to
\begin{align*}
\big(g_{i+1}+(q-q^{-1})\frac{ce_{i+1}}{b-c}\big)\big(g_{i}+(q-q^{-1})\frac{ce_{i}}{a-c}\big)\big(g_{i+1}+(q-q^{-1})\frac{be_{i+1}}{a-b}\big).
\end{align*}
By expanding the expression above, we get that it is equal to
\begin{align}\label{equation2}
g_{i+1}g_{i}g_{i+1}+&(q-q^{-1})\frac{bg_{i+1}g_{i}e_{i+1}}{a-b}+(q-q^{-1})\frac{cg_{i+1}e_{i}g_{i+1}}{a-c}\notag\\
+(q&-q^{-1})^{2}\frac{bcg_{i+1}e_{i}e_{i+1}}{(a-c)(a-b)}+(q-q^{-1})\frac{ce_{i+1}g_{i}g_{i+1}}{b-c}+
(q-q^{-1})^{2}\frac{bce_{i+1}g_{i}e_{i+1}}{(b-c)(a-b)}\notag\\
&+(q-q^{-1})^{2}\frac{c^{2}e_{i+1}e_{i}g_{i+1}}{(b-c)(a-c)}
+(q-q^{-1})^{3}\frac{bc^{2}e_{i+1}e_{i}e_{i+1}}{(b-c)(a-c)(a-b)}.
\end{align}

By \eqref{relations}, it is easy to see that \eqref{equation1} is equal to \eqref{equation2}. Thus, \eqref{Baxter-element1} holds.
\end{proof}

\section{Fusion formula for $Y_{r,n}$}

Largely inspired by [PA1, Section 7], we establish the fusion formula for the idempotent $E_{\mathcal{T}}$ in this section.

\subsection{Idempotents of $Y_{r,n}$}

We first recall the construction of the primitive idempotents $E_{\mathcal{T}}$ of $Y_{r, n}$ following [ChPA1], and then introduce the inductive formulae for $E_{\mathcal{T}}.$

For each $\bm{\nu}\in \mathcal{P}_{r,n}$, let $S^{\bm{\nu}}$ be the Specht module of $Y_{r,n}$ corresponding to $\bm{\nu}.$ Then the set $\{S^{\bm{\nu}}\:|\:\bm{\nu}\in \mathcal{P}_{r,n}\}$ forms a complete set of pairwise non-isomorphic irreducible representations of the split semisimple algebra $Y_{r, n}$. Assume that $m_{\bm{\nu}}=\dim S^{\bm{\nu}},$ and $S^{\bm{\nu}}$ has a basis $\{\bm{\text{v}}_{\mathcal{T}}\}$ indexed by $\mathrm{Std}(\bm{\nu})$. Thus, we can identify $\text{End}_{\mathbb{K}}(S^{\bm{\nu}})$ with the matrix algebra $\text{Mat}_{m_{\bm{\nu}}}(\mathbb{K})$, where each row and column of one matrix are labelled by the set $\{\bm{\text{v}}_{\mathcal{T}}\}.$ Since $Y_{r,n}$ is split semisimple over $\mathbb{K}$, by the Wedderburn-Artin theorem we get an isomorphism:
\begin{equation}\label{isomorphism}
I: Y_{r,n}\rightarrow \prod_{\bm{\nu}\in \mathcal{P}_{r,n}} \text{Mat}_{m_{\bm{\nu}}}(\mathbb{K}).
\end{equation}

For each $\bm{\nu}\in \mathcal{P}_{r,n}$, let $I_{\bm{\nu}}$ denote the projection $I_{\bm{\nu}}: Y_{r,n}\twoheadrightarrow \text{Mat}_{m_{\bm{\nu}}}(\mathbb{K}).$ For each $\bm{\lambda}\in \mathcal{P}_{r,n}$ and each $\mathcal{T}\in \mathrm{Std}(\bm{\lambda}),$ we continue to use the notations in \eqref{ci-pi-12345} in this section. Since $I$ is an isomorphism, we get a unique element $E_{\mathcal{T}}$ of $Y_{r,n}$ such that
\begin{align*}
I_{\bm{\nu}}(E_{\mathcal{T}})=
\begin{cases}
~~0 & \text{if } \bm{\nu}\neq\bm{\lambda};
\\
P_{\bm{\text{v}}_{\mathcal{T}}} & \text{otherwise,}
\end{cases}
\end{align*}
where $P_{\bm{\text{v}}_{\mathcal{T}}}$ is the diagonal $m_{\bm{\lambda}}\times m_{\bm{\lambda}}$ matrix with coefficient $1$ in the row labelled by $\bm{\text{v}}_{\mathcal{T}}$, and $0$ elsewhere. It follows from the definition of $E_{\mathcal{T}}$ and [ChPA1, Proposition 6] that we have
\begin{equation}\label{j-k-e-t}
J_{k}E_{\mathcal{T}}=E_{\mathcal{T}}J_{k}=\text{c}_{k}E_{\mathcal{T}}\quad\mbox{for}~k=1,\ldots,n.
\end{equation}
Moreover, it follows from [ChPA1, (5.2)] or [ChPA1, (7.5)-(7.6)] that we have
\begin{equation}\label{t-k-e-t}
t_{k}E_{\mathcal{T}}=E_{\mathcal{T}}t_{k}=\zeta_{\text{p}_{k}}E_{\mathcal{T}}\quad\mbox{for}~k=1,\ldots,n.
\end{equation}

Denote by $\bm{\theta}$ the $r$-node of $\mathcal{T}$ containing the number $n.$ Since the $r$-tableau $\mathcal{T}$ is standard, the $r$-node $\bm{\theta}$ is removable. Let $\mathcal{U}$ be the standard $r$-tableau obtained from $\mathcal{T}$ by removing $\bm{\theta}$ and let $\bm{\mu}$ be the shape of $\mathcal{U}$.

Denote by $\mathcal{E}_{+}(\bm{\mu})$ the set of $r$-nodes addable to $\bm{\mu}$. We have the following inductive formula for $E_{\mathcal{T}}$ in terms of the elements $t_1,\ldots,t_{n},J_1,\ldots,J_{n}$:
\begin{equation}\label{inductive-formula}
E_{\mathcal{T}}=E_{\mathcal{U}}\prod_{\substack{\bm{\theta'}\in \mathcal{E}_{+}(\bm{\mu})\\\text{c}(\bm{\theta'})\neq\text{c}(\bm{\theta})}}
\frac{J_{n}-\text{c}(\bm{\theta'})}{\text{c}(\bm{\theta})-\text{c}(\bm{\theta'})}\prod_{\substack{\bm{\theta'}\in \mathcal{E}_{+}(\bm{\mu})\\\text{p}(\bm{\theta'})\neq\text{p}(\bm{\theta})}}
\frac{t_{n}-\zeta_{\text{p}(\bm{\theta'})}}{\zeta_{\text{p}(\bm{\theta})}-\zeta_{\text{p}(\bm{\theta'})}}
\end{equation}
with $E_{\mathcal{T}_{0}}=1$ for the unique standard $r$-tableau $\mathcal{T}_{0}$ of size $0.$

\subsection{Fusion formulae of $E_{\mathcal{T}}$}

Recall that $S=\{\zeta_1,\ldots,\zeta_r\}$ is the set of all $r$-th roots of unity. We set
\begin{equation}\label{gamma-function}
\Gamma(v_1,\ldots,v_n) :=\prod_{i=1}^{n}\Big(\frac{\Pi_{\xi\in S}(v_i-\xi)}{v_i-t_i}\Big).
\end{equation}

Let $\phi_{1}(u) :=1$ and, for $k=2,\ldots,n$, set
\begin{align}\label{phi-function}
\phi_k(u_1,\ldots,u_{k-1},u)& :=g_{k-1}(u,u_{k-1})\phi_{k-1}(u_1,\ldots,u_{k-2},u)g_{k-1}^{-1}\notag\\
&=g_{k-1}(u,u_{k-1})g_{k-2}(u,u_{k-2})\cdots g_{1}(u,u_{1})\cdot g_{1}^{-1}\cdots g_{k-1}^{-1}.
\end{align}

We then define the following rational function:
\begin{align}\label{Phi-function}
\Phi(u_1,\ldots,u_n,v_1,\ldots,v_n) :=&\phi_n(u_1,\ldots,u_{n})\phi_{n-1}(u_1,\ldots,u_{n-1})\notag\\
&\cdots \phi_1(u_1)\Gamma(v_1,\ldots,v_n).
\end{align}

Now we can state the main result of this paper.

\begin{theorem}\label{main-theorem}
The idempotent $E_{\mathcal{T}}$ of $Y_{r,n}$ corresponding to the standard $r$-tableau $\mathcal{T}$ can be obtained by the following consecutive evaluations$:$
\begin{equation}\label{idempotents}
E_{\mathcal{T}}=\frac{1}{\emph{F}_{\bm{\lambda}}^{T}\emph{F}_{\bm{\lambda}}}\Phi(u_1,\ldots,u_{n},v_1,\ldots,v_{n})
\Big|_{v_{1}=\zeta_{\emph{p}_{1}}}\cdots\Big|_{v_{n}=\zeta_{\emph{p}_{n}}}\Big|_{u_{1}=\emph{c}_1}\cdots\Big|_{u_{n}=\emph{c}_{n}}.
\end{equation}
\end{theorem}

We first prove two necessary lemmas. We define the following element:
\begin{equation}\label{e-upn}
E_{\mathcal{U}, \text{p}_{n}} :=\frac{v-\zeta_{\text{p}_n}}{v-t_{n}}E_{\mathcal{U}}\Big|_{v=\zeta_{\text{p}_{n}}}.
\end{equation}
By definition, the element $E_{\mathcal{U}, \text{p}_{n}}$ is an idempotent which is equal to the sum of the idempotents $E_{\mathcal{S}}$, where $\mathcal{S}$ runs through the set of standard $r$-tableaux obtained from $\mathcal{U}$ by adding an $r$-node $\bm{\theta}$ with the number $n$ in it such that $\text{p}(\bm{\theta})=\text{p}_{n}.$

\begin{lemma}
Assume that $n\geq 1.$ We have
\begin{align}\label{F-PhiEu}
\emph{F}_{\mathcal{T}}(u)\phi_{n}(\emph{c}_1,\ldots,\emph{c}_{n-1},u)E_{\mathcal{U}, \emph{p}_{n}}=\frac{u-\emph{c}_n}{u-J_{n}}E_{\mathcal{U}, \emph{p}_{n}}.
\end{align}
\end{lemma}
\begin{proof}
We prove the lemma by induction on $n.$

When $n=1,$ the left-hand side of \eqref{F-PhiEu} is equal to $\frac{u-\text{c}_1}{u-1}\phi_{1}(u)E_{\mathcal{U}, \text{p}_{n}}=\frac{u-\text{c}_1}{u-1}E_{\mathcal{U}, \text{p}_{n}}$ by noting that $J_1=\phi_{1}(u)=1.$

When $n> 1$, we note that the left-hand side of \eqref{F-PhiEu} is equal to
\[\text{F}_{\mathcal{T}}(u)\Big(g_{n-1}+(q-q^{-1})\frac{\text{c}_{n-1}e_{n-1}}{u-{\text{c}_{n-1}}}\Big)\cdots \Big(g_{1}+(q-q^{-1})\frac{\text{c}_{1}e_{1}}{u-{\text{c}_{1}}}\Big)\cdot g_{1}^{-1}\cdots g_{n-1}^{-1}E_{\mathcal{U}, \text{p}_{n}}.\]
For $k=1,\ldots,n-1$, we have
\begin{align*}
e_{k}&\Big(g_{k-1}+(q-q^{-1})\frac{\text{c}_{k-1}e_{k-1}}{u-{\text{c}_{k-1}}}\Big)\cdots
\Big(g_{1}+(q-q^{-1})\frac{\text{c}_{1}e_{1}}{u-{\text{c}_{1}}}\Big)\\
&=\Big(g_{k-1}+(q-q^{-1})\frac{\text{c}_{k-1}e_{k-1}}{u-{\text{c}_{k-1}}}\Big)\cdots
\Big(g_{1}+(q-q^{-1})\frac{\text{c}_{1}e_{1}}{u-{\text{c}_{1}}}\Big)\cdot e_{1,k+1},
\end{align*}
and $e_{1,k+1}\cdot g_{1}^{-1}\cdots g_{n-1}^{-1}=g_{1}^{-1}\cdots g_{n-1}^{-1}\cdot e_{k,n}.$

Moreover, $e_{k,n} E_{\mathcal{U}, \text{p}_{n}}=0$ if $\text{p}_{k}\neq\text{p}_{n}$. Thus the left-hand side of \eqref{F-PhiEu} is equal to
\begin{align}\label{Fg-PhiEu}
\text{F}_{\mathcal{T}}(u)&\Big(g_{n-1}+(q-q^{-1})\frac{\delta_{\text{p}_{n-1}, \text{p}_{n}}\text{c}_{n-1}e_{n-1}}{u-{\text{c}_{n-1}}}\Big)\cdots\notag\\ &\times\Big(g_{1}+(q-q^{-1})\frac{\delta_{\text{p}_{1}, \text{p}_{n}}\text{c}_{1}e_{1}}{u-{\text{c}_{1}}}\Big)\cdot g_{1}^{-1}\cdots g_{n-1}^{-1}E_{\mathcal{U}, \text{p}_{n}}.
\end{align}

Suppose first that $\text{p}_i\neq \text{p}_n$ for $i=1,\ldots,n-1.$ In this situation, we have $E_{\mathcal{U}, \text{p}_{n}}=E_{\mathcal{T}}$ and $\text{c}_n=1.$ Thus, we have $\text{F}_{\mathcal{T}}(u)=\frac{u-\text{c}_{n}}{u-1}=1.$ Due to \eqref{Fg-PhiEu}, we have the left-hand side of \eqref{F-PhiEu} is equal to $\text{F}_{\mathcal{T}}(u)E_{\mathcal{U}, \text{p}_{n}}=E_{\mathcal{T}}$; while the right-hand side of \eqref{F-PhiEu} is also equal to $\frac{u-\text{c}_{n}}{u-\text{c}_{n}}E_{\mathcal{U}, \text{p}_{n}}=E_{\mathcal{T}}$.

Next assume that there exists some $l\in \{1,\ldots,n-1\}$ such that $\text{p}_{l}=\text{p}_{n}.$ Fix $l$ such that $\text{p}_{l}=\text{p}_{n}$ and $\text{p}_{i}\neq \text{p}_{n}$ for $i=l+1,\ldots,n-1.$

Let $\mathcal{V}$ be the standard $r$-tableau obtained from $\mathcal{U}$ by removing the $r$-nodes with numbers $l+1,\ldots,n-1$ and $\mathcal{W}$ be the standard $r$-tableau obtained from $\mathcal{V}$ by removing the $r$-node with the number $l.$ We define
\begin{equation*}
E_{\mathcal{W}, \text{p}_{l}} :=\frac{v-\zeta_{\text{p}_l}}{v-t_{l}}E_{\mathcal{W}}\Big|_{v=\zeta_{\text{p}_{l}}}.
\end{equation*}

Since $E_{\mathcal{W}}$ can be expressed in terms of $J_1,\ldots,J_{l-1}$ and $t_1,\ldots,t_{l-1}$, $E_{\mathcal{W}}$ commutes with $g_{l}^{-1}g_{l+1}^{-1}\cdots g_{n-1}^{-1}$. Note that $E_{\mathcal{W}}E_{\mathcal{U}}=E_{\mathcal{U}}=E_{\mathcal{U}}^{2},$ $E_{\mathcal{U}, \text{p}_{n}}^{2}=E_{\mathcal{U}, \text{p}_{n}}$, $\text{p}_{l}=\text{p}_{n}$ and $t_lg_{l}^{-1}\cdots g_{n-1}^{-1}$$=g_{l}^{-1}\cdots g_{n-1}^{-1}t_n.$ Thus, we have
\begin{align}\label{EW-PhiEu}
E_{\mathcal{W}, \text{p}_{l}}g_{l}^{-1}g_{l+1}^{-1}\cdots g_{n-1}^{-1}E_{\mathcal{U}, \text{p}_{n}}&=g_{l}^{-1}\cdots g_{n-1}^{-1}\frac{v-\zeta_{\text{p}_l}}{v-t_{n}}E_{\mathcal{W}}E_{\mathcal{U}, \text{p}_{n}}\Big|_{v=\zeta_{\text{p}_{l}}}\notag\\
&=g_{l}^{-1}\cdots g_{n-1}^{-1}E_{\mathcal{U}, \text{p}_{n}}.
\end{align}
By \eqref{EW-PhiEu}, we rewrite \eqref{Fg-PhiEu} as follows:
\begin{align*}
\text{F}_{\mathcal{T}}(u)g_{n-1}\cdots g_{l+1}\Big(g_{l}+(q-q^{-1})\frac{\text{c}_{l}e_{l}}{u-{\text{c}_{l}}}\Big)\phi_{l}(\emph{c}_1,\ldots,\emph{c}_{l-1},u)E_{\mathcal{W}, \text{p}_{l}}g_{l}^{-1}\cdots g_{n-1}^{-1}E_{\mathcal{U}, \emph{p}_{n}}.
\end{align*}
By the induction hypothesis, we have
\begin{align*}
\phi_{l}(\text{c}_1,\ldots,\text{c}_{l-1},u)E_{\mathcal{W}, \text{p}_{l}}=\text{F}_{\mathcal{V}}(u)^{-1}\frac{u-\text{c}_l}{u-J_{l}}E_{\mathcal{W}, \text{p}_{l}},
\end{align*}
and we use \eqref{EW-PhiEu} again to obtain that the left-hand side of \eqref{F-PhiEu} is equal to
\begin{align}\label{FF-PhiEu}
\text{F}_{\mathcal{T}}(u)\text{F}_{\mathcal{V}}(u)^{-1}g_{n-1}\cdots g_{l+1}\Big(g_{l}+(q-q^{-1})\frac{\text{c}_{l}e_{l}}{u-{\text{c}_{l}}}\Big)\frac{u-\text{c}_l}{u-J_{l}}g_{l}^{-1}\cdots g_{n-1}^{-1}E_{\mathcal{U}, \emph{p}_{n}}.
\end{align}

Noting that $J_{n}$ commutes with $E_{\mathcal{U}, \emph{p}_{n}},$ we can move $(u-J_n)^{-1}$ from the right-hand side of \eqref{F-PhiEu} to the left-hand side. Since $g_{k}^{-1}=g_k-(q-q^{-1})e_{k}$, $e_kg_{k+1}\cdots g_{n-1}=g_{k+1}\cdots g_{n-1}e_{k,n}$ and $e_{k,n} E_{\mathcal{U}, \text{p}_{n}}=0$ for $k=l+1,\ldots,n-1,$ we can move $g_{n-1}\cdots g_{l+1}$ to the right-hand side. By \eqref{Baxter-element2}, $g_{l}(u, \text{c}_l)$ is invertible. We finally get that \eqref{F-PhiEu} is equivalent to
\begin{align}\label{Equi-equality}
\text{F}_{\mathcal{T}}(u)&\text{F}_{\mathcal{V}}(u)^{-1}(u-\text{c}_l)g_{l}^{-1}\cdots g_{n-1}^{-1}(u-J_n)E_{\mathcal{U}, \text{p}_{n}}
=(u-\text{c}_n)(u-J_l)\notag\\
&\times\Big(g_{l}+(q-q^{-1})\frac{ue_{l}}{\text{c}_{l}-u}\Big)
\Big(1-(q-q^{-1})^{2}\frac{u\text{c}_le_{l}}{(u-\text{c}_{l})^{2}}\Big)^{-1}g_{l+1}\cdots g_{n-1}E_{\mathcal{U}, \text{p}_{n}}.
\end{align}

Since $\text{p}_{l}=\text{p}_{n}$ and $\text{p}_{i}\neq \text{p}_{n}$ for $i=l+1,\ldots,n-1,$ we have, by the definition \eqref{f-T}, that
\begin{align*}
\text{F}_{\mathcal{T}}(u)\text{F}_{\mathcal{V}}(u)^{-1}=\frac{u-\text{c}_n}{u-\text{c}_l}\frac{(u-\text{c}_l)^{2}}{(u-\text{c}_l)^{2}-(q-q^{-1})^{2}u\text{c}_l}.
\end{align*}
Notice that $e_lg_{l+1}\cdots g_{n-1}=g_{l+1}\cdots g_{n-1}e_{l,n}$ and that $e_{l,n}E_{\mathcal{U}, \text{p}_{n}}=E_{\mathcal{U}, \text{p}_{n}}$ since $\text{p}_{l}=\text{p}_{n}.$ Therefore, to verify \eqref{Equi-equality}, it suffices to show that
\begin{align}\label{Equi-equality1}
g_{l}^{-1}\cdots g_{n-1}^{-1}(u-J_n)E_{\mathcal{U}, \text{p}_{n}}
=(u-J_l)\Big(g_{l}+(q-q^{-1})\frac{ue_{l}}{\text{c}_{l}-u}\Big)g_{l+1}\cdots g_{n-1}E_{\mathcal{U}, \text{p}_{n}}.
\end{align}

By \eqref{JM-elements}, we have $g_{l}^{-1}g_{l+1}^{-1}\cdots g_{n-1}^{-1}J_n=J_lg_{l}g_{l+1}\cdots g_{n-1}$. Thus, the left-hand side of \eqref{Equi-equality1} is equal to
\begin{align}\label{Equi-equality2}
ug_{l}^{-1}\cdots g_{n-1}^{-1}E_{\mathcal{U}, \text{p}_{n}}-J_lg_{l}\cdots g_{n-1}E_{\mathcal{U}, \text{p}_{n}}.
\end{align}
By the fact that $e_{l,n} E_{\mathcal{U}, \text{p}_{n}}=E_{\mathcal{U}, \text{p}_{n}}$ and $e_{k,n} E_{\mathcal{U}, \text{p}_{n}}=0$ for $k=l+1,\ldots,n-1,$ we get that \eqref{Equi-equality2} is equal to
\begin{align}\label{Equi-equality3}
&ug_{l}^{-1}g_{l+1}\cdots g_{n-1}E_{\mathcal{U}, \text{p}_{n}}-J_lg_{l}\cdots g_{n-1}E_{\mathcal{U}, \text{p}_{n}}\notag\\
=&u\big(g_l-(q-q^{-1})e_l\big)g_{l+1}\cdots g_{n-1}E_{\mathcal{U}, \text{p}_{n}}-J_lg_{l}\cdots g_{n-1}E_{\mathcal{U}, \text{p}_{n}}\notag\\
=&(u-J_l)g_{l}\cdots g_{n-1}E_{\mathcal{U}, \text{p}_{n}}-(q-q^{-1})ug_{l+1}\cdots g_{n-1}E_{\mathcal{U}, \text{p}_{n}}.
\end{align}

By \eqref{giXj} we see that $J_l$ commutes with $g_{l+1}\cdots g_{n-1}$. Thus, by $J_lE_{\mathcal{U}, \text{p}_{n}}=\text{c}_lE_{\mathcal{U}, \text{p}_{n}},$ we see that the right-hand side of \eqref{Equi-equality1} is equal to
\begin{align}\label{Equi-equality4}
&(u-J_l)g_{l}\cdots g_{n-1}E_{\mathcal{U}, \text{p}_{n}}+(q-q^{-1})(u-J_l)\frac{ue_{l}}{\text{c}_{l}-u}g_{l+1}\cdots g_{n-1}E_{\mathcal{U}, \text{p}_{n}}\notag\\
=&(u-J_l)g_{l}\cdots g_{n-1}E_{\mathcal{U}, \text{p}_{n}}+(q-q^{-1})(u-J_l)\frac{u}{\text{c}_{l}-u}g_{l+1}\cdots g_{n-1}E_{\mathcal{U}, \text{p}_{n}}\notag\\
=&(u-J_l)g_{l}\cdots g_{n-1}E_{\mathcal{U}, \text{p}_{n}}+(q-q^{-1})\cdot\frac{u}{\text{c}_{l}-u}g_{l+1}\cdots g_{n-1}(u-J_l)E_{\mathcal{U}, \text{p}_{n}}\notag\\
=&(u-J_l)g_{l}\cdots g_{n-1}E_{\mathcal{U}, \text{p}_{n}}+(q-q^{-1})\cdot\frac{u}{\text{c}_{l}-u}g_{l+1}\cdots g_{n-1}(u-\text{c}_l)E_{\mathcal{U}, \text{p}_{n}}\notag\\
=&(u-J_l)g_{l}\cdots g_{n-1}E_{\mathcal{U}, \text{p}_{n}}-(q-q^{-1})ug_{l+1}\cdots g_{n-1}E_{\mathcal{U}, \text{p}_{n}}.
\end{align}

Comparing \eqref{Equi-equality4} with \eqref{Equi-equality3}, we see that \eqref{Equi-equality1} holds.
\end{proof}

Recall the definition of $\phi_k(u_1,\ldots,u_{k-1},u)$ in \eqref{phi-function}. For $k=1,\ldots,n,$ we define
\begin{equation}\label{tilde-phi}
\widetilde{\phi}_{k}(u_1,\ldots,u_{k-1},u,v) :=\phi_k(u_1,\ldots,u_{k-1},u)\cdot\Big(\frac{\Pi_{\xi\in S}(v-\xi)}{v-t_k}\Big).
\end{equation}

\begin{lemma}
Assume that $n\geq 1.$ We have
\begin{align}\label{FF-Phi}
\emph{F}_{\mathcal{T}}^{T}(v)\emph{F}_{\mathcal{T}}(u)\widetilde{\phi}_{n}(\emph{c}_1,\ldots,\emph{c}_{n-1},u,v)
E_{\mathcal{U}}\Big|_{v=\zeta_{\emph{p}_{n}}}=\frac{u-\emph{c}_n}{u-J_{n}}\frac{v-\zeta_{\emph{p}_n}}{v-t_{n}}E_{\mathcal{U}}\Big|_{v=\zeta_{\emph{p}_{n}}}.
\end{align}
\end{lemma}
\begin{proof}
By \eqref{f-TT}, we have
\begin{align}\label{Ftt-equality}
\emph{F}_{\mathcal{T}}^{T}(v)\cdot\Big(\frac{\Pi_{\xi\in S}(v-\xi)}{v-t_n}\Big)=\frac{v-\zeta_{\text{p}_n}}{v-t_{n}}.
\end{align}
By \eqref{e-upn}, \eqref{tilde-phi} and \eqref{Ftt-equality}, we see that \eqref{FF-Phi} is a direct consequence of \eqref{F-PhiEu}.
\end{proof}

{\it Proof of Theorem \ref{main-theorem}}
Since $g_i$ commutes with $t_k$ if $i< k-1$, we can rewrite the function $\Phi(u_1,\ldots,u_n,v_1,$$\ldots,v_n)$ as follows:
\begin{align}\label{Phitilde}
\Phi(&u_1,\ldots,u_n,v_1,\ldots,v_n)\notag\\
&=\widetilde{\phi}_{n}(u_1,\ldots,u_{n},v_{n})\widetilde{\phi}_{n-1}(u_1,\ldots,u_{n-1},v_{n-1})\cdots \widetilde{\phi}_{1}(u_1,v_{1}).
\end{align}

We prove this theorem by induction on $n.$ For $n=0,$ the situation is trivial.

For $n> 0,$ by \eqref{Phitilde} and the induction hypothesis we can rewrite the right-hand side of \eqref{idempotents} as follows:
\[(\text{F}_{\bm{\lambda}}^{T}\text{F}_{\bm{\lambda}})^{-1}\text{F}_{\bm{\mu}}^{T}\text{F}_{\bm{\mu}}
\widetilde{\phi}_{n}(\text{c}_1,\ldots,\text{c}_{n-1},u_n,v_n)E_{\mathcal{U}}\Big|_{v_n=\zeta_{\text{p}_{n}}}\Big|_{u_n=\text{c}_{n}}.\]
By \eqref{FF-Phi} we can rewrite the expression above as
\begin{equation}\label{add-equation}
(\text{F}_{\bm{\lambda}}^{T}\text{F}_{\bm{\lambda}})^{-1}\text{F}_{\bm{\mu}}^{T}\text{F}_{\bm{\mu}}
(\text{F}_{\mathcal{T}}^{T}(v_n)\text{F}_{\mathcal{T}}(u_n))^{-1}
\frac{u_n-\text{c}_n}{u_n-J_{n}}\frac{v_n-\zeta_{\text{p}_n}}{v_n-t_{n}}E_{\mathcal{U}}\Big|_{v_n=\zeta_{\text{p}_{n}}}\Big|_{u_n=\text{c}_{n}}.
\end{equation}

Assume that $\{\mathcal{T}_1,\ldots,\mathcal{T}_k\}$ is the set of pairwise different standard $r$-tableaux obtained from $\mathcal{U}$ by adding an $r$-node containing the number $n.$ Notice that $\mathcal{T}\in \{\mathcal{T}_1,\ldots,\mathcal{T}_k\}.$ Moreover, we have
\begin{equation}\label{sum-formula}
E_{\mathcal{U}}=\sum_{i=1}^{k}E_{\mathcal{T}_{i}}.
\end{equation}
Consider the following rational function in $u$ and $v$:
\begin{equation}\label{rational-function}
\frac{u-\text{c}_n}{u-J_{n}}\frac{v-\zeta_{\text{p}_n}}{v-t_{n}}E_{\mathcal{U}}.
\end{equation}
The formulae \eqref{j-k-e-t} and \eqref{t-k-e-t} imply that \eqref{rational-function} is non-singular at $u=\text{c}_n$ and $v=\zeta_{\text{p}_n}$, and moreover, by replacing $E_{\mathcal{U}}$ with the right-hand side of \eqref{sum-formula}, we get
\begin{equation}\label{sum-function}
\frac{u-\text{c}_n}{u-J_{n}}\frac{v-\zeta_{\text{p}_n}}{v-t_{n}}E_{\mathcal{U}}\Big|_{v=\zeta_{\text{p}_{n}}}\Big|_{u=\text{c}_n}=E_{\mathcal{T}}.
\end{equation}

By \eqref{f-T-relation} and \eqref{f-T-u}, together with \eqref{add-equation} and \eqref{sum-function}, we see that the right-hand side of \eqref{idempotents} is equal to $E_{\mathcal{T}}.$ $\hfill{} \Box$

Finally, let us look at an example.
\begin{example}
Consider the situation that $r=n=3$ and the $3$-partition $\bm{\lambda}=((2),(0),(1))$ of $3$. We shall consider the following standard $3$-tableau of shape $\bm{\lambda}$:
\[\mathcal{T}=\left(\hspace{-0.1cm}\begin{array}{l}\fbox{1}\fbox{3}\\[-0.05em] \end{array}
\,,\,\varnothing\,,\, \begin{array}{l}\fbox{2}\\[-0.05em] \end{array}\hspace{-0.1cm}\right).\]
Note that
\begin{align}\label{v-zetavt}
\frac{\Pi_{\zeta\in S}(v_i-\zeta)}{v_i-t_i}=\frac{v_{i}^{r}-1}{v_i-t_i}=\frac{v_{i}^{r}-t_{i}^{r}}{v_i-t_i}=\sum_{s=0}^{r-1}v_{i}^{s}t_{i}^{r-1-s}.
\end{align}
Theorem \ref{main-theorem} implies that the idempotent $E_{\mathcal{T}}$ can be expressed as
\begin{align*}
E_{\mathcal{T}}=&\frac{q\zeta_{1}^{2}\zeta_{3}}{27(q+q^{-1})}g_{2}(q^{2}, 1)g_{1}(q^{2}, 1)g_{1}^{-1}g_{2}^{-1}\times g_{1}(1,1)g_{1}^{-1}\\
&\times(\zeta_{1}^{2}+\zeta_{1}t_{1}+t_{1}^{2})(\zeta_{3}^{2}+\zeta_{3}t_{2}+t_{2}^{2})(\zeta_{1}^{2}+\zeta_{1}t_{3}+t_{3}^{2}).
\end{align*}
\end{example}

\noindent{\bf Acknowledgements.}
The author is deeply indebted to Dr. Shoumin Liu for posing this question to him and for helpful comments. Many ideas of this paper originate from the reference [PA1].
%/////////////////////////////////////////////////////////////////////////////////////////////////////////////////////////////////////////////////////////

%by (\ref{stn})
%\begin{align*}
%  t_{n+1}&\geq\frac{1}{t_n+(1-t_n)\frac{1}{c}}\\
%  1-t_{n+1}&\leq \frac{(1-t_n)(\frac{1}{c}-1)}{t_n+(1-t_n)\frac{1}{c}}=\frac{(1-t_n)(\frac{1}{c}-1)}{(1-t_n)(\frac{1}{c}-1)+1}\\
%  \frac{1}{1-t_{n+1}}&\geq \frac{1}{(1-t_n)(\frac{1}{c}-1)+1}
%\end{align*}
%Since $c\geq \frac{1}{2}$, $\frac{1}{c}-1\geq 1$, we have
%$$\frac{1}{1-t_{n+1}}\geq \frac{1}{1-t_n}+1\geq\cdots\geq\frac{1}{1-t_1}+n\geq n+1$$
%i.e. $1-t_{n+1}\leq \frac{1}{n+1}$
%$$0\leq v_{2n+1}-x^*\leq v_{2n+1}-v_{2n}\leq v_{2n+1}-t_nv_{2n+1}\geq \frac{1}{n}v_{2n+1}$$
%So $\|v_{2n+1}-x^*\|\leq \frac{N}{n}\|\bar{v}\|$, where $N$ is the normal constant of $P$.
%/////////////////////////////////////////////////////////////////////////////////////////////////////////////////////////////////

\end{document}